\documentclass[12pt]{article}
\usepackage{amsxtra,amssymb,amsthm,amsmath,latexsym}

\theoremstyle{plain}
\newtheorem{theorem}{Theorem}

\newtheorem{lemma}{Lemma}
\newtheorem*{A}{Assumption A}
\numberwithin{equation}{section}

\newcommand{\refS}[1]{Section~\ref{S:#1}}

\def\qed{{\hfill $\Box$}}

\def\R{{\mathbb R}}

\def\oH1{{\overset{\circ}{H}\kern-.02in{}^1}}

\def\const{{\,const\,}}

\def\bee{\begin{equation*}}\def\eee{\end{equation*}}
\def\be{\begin{equation}}\def\ee{\end{equation}}

\begin{document}
\title{Existence of a solution to a nonlinear equation}
\author{A.G. Ramm\\
 Mathematics Department, Kansas State University, \\
 Manhattan, KS 66506-2602, USA\\
ramm@math.ksu.edu\\}

\date{}
\maketitle\thispagestyle{empty}

\begin{abstract} \footnote{Math subject classification: 35J65 }
\footnote{key words: nonlinear elliptic equations } Equation
$(-\Delta+k^2)u+f(u)=0$ in $D$, $u\mid_{\partial D}=0$, where $k=\const>0$
and $D\subset\R^3$ is a bounded domain, has a solution if $f:\R\to\R$ is a
continuous function in the region $|u|\geq a$, piecewise-continuous
 in the region $|u|\leq a$, with finitely many discontinuity points $u_j$
such that $f(u_j\pm 0)$ exist, and 
 $uf(y)\geq 0$ for $|u|\geq a$, where $a\geq
0$ is an arbitrary fixed number. \end{abstract} 

\section{Introduction}\label{S:1} Let $D\subset\R^3$ be a bounded domain
with Lipschitz boundary $S$, $k=\const>0$ , $f:\R\to\R$ be a
function such that \be\label{e1.1} uf(u)\geq 0, \hbox{\quad
for\quad} |u|\geq a\geq 0,\ee 
where $a$ is an arbitrary fixed number, and $f$ is continuous in the 
region $|u|\geq a$, and bounded and piecewise-continuous with at most 
finitely many discontinuity points $u_j$, such that 
$f(u_j+0)$ and $f(u_j-0)$ exist, in the region $|u|\leq a$.

Consider the problem \be\label{e1.2} (-\Delta+k^2)u+f(u)=0 \hbox{\quad
in\quad} D, \ee \be\label{e1.3} u=0\hbox{\quad on\quad} S.\ee There is a
large literature on problems of this type. Usually it is assumed that $f$
does not grow too fast or $f$ is 
monotone (see e.g. \cite{B} and references therein). 

 {\it The novel point in
this note  is the absence of monotonicity restrictions 
on $f$ and of the growth restrictions on $f$ as $|u|\to \infty$, except 
for the assumption (1.1)}. 

This assumption allows an arbitrary behavior of 
$f$
inside the region $|u|\leq a$, where $a\geq 0$ can be arbitrary large,
and an arbitrarily rapid growth of $f$ to $+\infty$ as $u\to +\infty$, or 
arbitrarily rapid 
decay of  $f$ to $-\infty$ as $u\to -\infty$. 

Our result is:

\begin{theorem}\label{T:1}
Under the above assumptions problem \eqref{e1.2}--\eqref{e1.3} has 
a solution
$u\in H^2(D)\cap\oH1(D):=H^2_0(D)$.
\end{theorem}

Here $H^\ell(D)$ is the usual Sobolev space, $\oH1(D)$ is the closure of
$C^\infty_0(D)$ in the norm $H^1(D)$. Uniqueness of the solution does not
hold without extra assumptions.

The ideas of our proof are: first, we prove that if 
$\sup_{u\in \R}|f(u)|\leq\mu$, 
then a solution to \eqref{e1.2}--\eqref{e1.3} exists by the Schauder's 
fixed-point theorem. Here
$\mu$ is a constant. Secondly, we prove an a priori bound 
$\|u\|_\infty\leq a$.
If this bound is proved,
then the solution to problem \eqref{e1.2}--\eqref{e1.3} with $f$ replaced
by \be\label{e1.4} F(u):=\begin{cases}f(u), & |u|\leq a\\f(a),&u\geq
a\\f(-a),&u\leq -a\end{cases}\ee has a solution, and this solution solves
the original problem \eqref{e1.2}--\eqref{e1.3}. The bound 
$\|u\|_\infty\leq a$ is proved by using some integral inequalities.
An alternative proof of this bound is also given. This proof is based on 
the maximum 
principle for elliptic equation (1.2).

In \refS{2} proofs are given.
We use some ideas from \cite{R129}.

\section{Proofs.}\label{S:2} 
If $u\in L^\infty:=L^\infty(D)$, then problem
\eqref{e1.2}--\eqref{e1.3} is equivalent to the integral equation:
\be\label{e2.1}u=-\int_D G(x,y)f(u(y))dy:=T(u),\ee where \be\label{e2.2}
(-\Delta+k^2)G=-\delta(x-y)\hbox{\quad in\quad} D,\qquad g\mid_{x\in S}
=0.\ee By the maximum principle, \be\label{e2.3} 0\leq
G(x,y)<g(x,y):=\frac{e^{-k|x-y|}}{4 \pi|x-y|},\qquad x,y\in D.\ee The map
$T$ is a continuous and compact map in the space $C(D):=X$, and 
\be\label{e2.4} \|u\|_{C(D)}:= \|u\|\leq \mu \sup_x \int_D
\frac{e^{-k|x-y|}}{4 \pi|x-y|}dy \leq \mu \int_{\R^3} \frac{e^{-k|y|}}{4
\pi|y|} dy\leq\frac{\mu}{k^2}\ee This is an a priori estimate of any
bounded solution to \eqref{e1.2}--\eqref{e1.3} for a bounded nonlinearity
$f$ such that $\sup_{u\in \R}|f(u)|\leq \mu$. Thus, Schauder's fixed-point 
theorem yields the existence
of a solution to \eqref{e2.1}, and consequently to problem
\eqref{e1.2}--\eqref{e1.3}, for bounded $f$. Indeed, if $B$ is a closed
ball of radius $\frac{\mu}{k^2}$, then the map $T$ maps this ball into
itself by \eqref{e2.4}, and since $T$ is compact, the Schauder principle
is applicable. Thus, the following lemma is proved.

\begin{lemma}\label{L:1}  
If $\sup_{u\in\R}|f(u)|\leq\mu$, then problems \eqref{e2.1} and 
\eqref{e1.2}--\eqref{e1.3}
have a solution in $C(D)$, and this solution satisfies estimate \eqref{e2.4}.
\end{lemma}

Let us now prove an a priori bound for any solution $u\in C(D)$ of the
problem \eqref{e1.2}--\eqref{e1.3} without assuming that
$\sup_{u\in \R}|f(u)|<\infty$.

Let $u_+:=\max(u,0)$, $u_-=\max(-u,0)$. Multiply \eqref{e1.2} by 
$(u-a)_+$, integrate over
$D$ and then by parts to get
\be\label{e2.5}
0=\int_D [\nabla u\cdot\nabla(u-a)_+ +k^2u(u-a)_+ +f(u)(u-a)_+]dx,\ee
where we have integrated by parts and the boundary integral vanishes 
because 
$(u-a)_+=0$ on $S$ for $a\geq 0$. Each of the terms in \eqref{e2.5} is 
nonnegative, the last one due to \eqref{e1.1}.
Thus \eqref{e2.5} implies
\be\label{e2.6} u\leq a.\ee
Similarly, using \eqref{e1.1} again, and multiplying \eqref{e1.2} by 
$(-u-a)_+$, one gets
\be\label{e2.7}-a\leq u.\ee

We have proved:

\begin{lemma}\label{L:2} If \eqref{e1.1} holds, then any solution $u\in
H^2_0(D)$ to \eqref{e1.2}--\eqref{e1.3} satisfies the inequality
\be\label{e2.8} |u(x)|\leq a.\ee \end{lemma}

Consider now equation \eqref{e2.1} in $C(D)$ with an arbitrary continuous
$f$ satisfying \eqref{e1.1}.  Any $u\in C(D)$ which solves \eqref{e2.1}
solves \eqref{e1.2}--\eqref{e1.3} and therefore satisfies \eqref{e2.8} and
belongs to $H^2_0(D)$. This $u$ solves  problem
\eqref{e1.2}--\eqref{e1.3} with $f$ replaced by $F$ defined in
\eqref{e1.4}, and vice-versa.  Since $F$ is a bounded nonlinearity,
equation \eqref{e2.1} and problem \eqref{e1.2}--\eqref{e1.3} (with $f$
replaced by $F$) has a solution by Lemma 1.

Theorem 1 is proved. \qed

{\it An alternative proof of the estimate (2.8):
}

Let us sketch an alternative derivation of the inequality (2.8) using the 
maximum principle. Let us derive (2.6). The derivation of (2.7) is 
similar. 

Assume that (2.6) fails. Then $u>a$ at some point in $D$. Therefore  
at a point $y$ at which $u$ attains its maximum value one has
$u(y)\geq u(x)$ for all $x\in D$  and 
$u(y)>a$. The function $u$ attains its maximum value, which is positive, 
at some point in $D$, because $u$ is continuous, vanishes at the boundary 
of $D$, and is positive at some point of $D$ by the assumption $u>a$.
At the point $y$, where the function $u$ 
attains its maximum,
one has $-\Delta u\geq 0$ and $k^2 u(y)>0$. Moreover, $f(u(y))>0$ by the 
assumption 
(1.1), since $u(y)>a$. Therefore the left-hand side of  
equation (1.2) is positive, while its left-hand side is zero.
Thus we have got a 
contradiction, and the estimate (2.6) is proved. Similarly one proves 
estimate (2.7). Thus, (2.8) is proved. \qed

\end{document}